\documentclass[12pt, dvips, twoside]{amsart}
\usepackage{amsmath,amssymb,amscd,amsxtra,latexsym,bbm,graphicx,epsfig,epic,
eepic,oldgerm,bm,euscript,psfrag,amsthm}
\usepackage{a4wide}
\usepackage[dvips]{color}

\input xy
\xyoption{all}
\CompileMatrices

\theoremstyle{plain}
\newtheorem{theorem}{Theorem}[section]
\newtheorem*{theorem*}{Theorem}
\newtheorem{lemma}[theorem]{Lemma}
\newtheorem{proposition}[theorem]{Proposition}

\newtheorem*{remark*}{Remark}
\newtheorem*{remarks*}{Remarks}
\newtheorem{remark}[theorem]{Remark}

\newtheorem{claim}[theorem]{Claim}

\newtheorem*{example*}{Example}
\newtheorem*{examples*}{Examples}

\newtheorem*{definition*}{Definition}

\newcommand{\proofend}{\hspace*{\fill} $\Box$\\}
\newcommand{\diam}{\hspace*{\fill} $\Diamond$\\}

\def\1{\:\!}
\def\2{\;\!}
\def\s{\smallskip}
\def\m{\medskip}
\def\Ker{\operatorname{Ker}}
\def\Im{\operatorname {Im}}

\def\Diffc0{\operatorname{Diff^c_0}}

\def\Sympc0{\operatorname{Symp^c_0}}

\def\const{\operatorname{const}}

\def\Ai{\operatorname{Ai}}
\def\codim{\operatorname{codim}}

\def\ga{\alpha}
\def\gb{\beta}

\def\eps{\varepsilon}

\def\gf{\varphi}

\def\gl{\lambda}

\def\gs{\sigma}

\def\cn{{\mathcal N}}

\def\NN{\mathbbm{N}}

\def\RR{\mathbbm{R}}

\def\ZZ{\mathbbm{Z}}

\def\pp{\partial}

\def\lef{\operatorname{Left}}
\def\righ{\operatorname{Right}}

\def\ni{\noindent}
\def\b{\bigskip}
\def\m{\medskip}

\def\.{\mskip1mu}
\def\?{\mskip-1mu}

\def\proof{\noindent {\it Proof. \;}}

\begin{document}

\title{Bifurcating extremal domains for the first eigenvalue of the Laplacian}

\author{Felix Schlenk}
\address{(F.~Schlenk) 
Institut de Math\'ematiques,
Universit\'e de Neuch\^atel, 
Rue \'Emile Argand~11, 
CP~158,
2009 Neuch\^atel,
Switzerland} 
\email{schlenk@unine.ch}

\author{Pieralberto Sicbaldi}
\address{(P.~Sicbaldi)
Laboratoire d'Analyse Topologie Probabilit\'es,
Universit\'e Aix-Marseille~3,
Avenue de l'Escadrille Normandie Niemen, 
13397 Marseille cedex 20,
France}
\email{pieralberto.sicbaldi@univ-cezanne.fr}

\date{\today}

\begin{abstract}
We prove the existence of a smooth family of non-compact domains 
$\Omega_s \subset \RR^{n+1}$, $n \ge 1$,
bifurcating from the straight cylinder $B^{n} \times \RR$ 
for which the first eigenfunction of the Laplacian with $0$~Dirichlet boundary condition 
also has constant Neumann data at the boundary:
For each~$s \in (-\eps,\eps)$, the overdetermined system
$$
\left\{
\begin{array} {ll}
\Delta\, u + \gl\, u = 0 &\mbox{in }\; \Omega_s\\
u=0 & \mbox{on }\; \pp \Omega_s \\
\langle \nabla u, \nu \rangle = \const &\mbox{on }\; \pp \Omega_s 
\end{array} 
\right.
$$
has a bounded positive solution. 
The domains $\Omega_s$ are rotationally symmetric and periodic 
with respect to the $\RR$-axis of the cylinder; 
they are of the form
$$
\Omega_s \,=\,
\left\{ 
(x,t) \in \RR^{n} \times \RR  \,\, \mid \,\, \|x\| < 
1+s \cos \left( \frac{2\pi}{T_s}\,t \right) + O(s^2)  
\right\} 
$$
where $T_s = T_0 + O(s)$ and $T_0$ is a positive real number depending on~$n$.
For $n \ge 2$
these domains provide a smooth family of counter-examples to a conjecture of 
Berestycki, Caffarelli and Nirenberg. 
We also give rather precise upper and lower bounds for the bifurcation period~$T_0$. 
This work improves a recent result of the second author.
\end{abstract}

\maketitle


\section{Introduction and main results}  \label{s:intro}

\subsection{The problem}  
\label{ss:problem}
Let $\Omega$ be a bounded domain in $\RR^n$ with smooth boundary,
and consider the Dirichlet problem
\begin{equation} \label{e:sys1}
\left\{\begin{array} {ll}
\Delta\, u + \gl\, u = 0 &\mbox{in }\; \Omega\\
u=0 & \mbox{on }\; \pp \Omega .
\end{array}\right.
\end{equation}
Denote by $\gl_1(\Omega)$ the smallest positive constant $\gl$ for which this system has a solution (i.e.~$\gl_1(\Omega)$ is the first eigenvalue of the Laplacian on $\Omega$ with 0 Dirichlet boundary condition).
By the Krein--Rutman theorem, 
the corresponding solution~$u$ 
(i.e.~the first eigenfunction of the Laplacian on $\Omega$ with 0~Dirichlet boundary condition) 
is positive on~$\Omega$, and $u$ is the only eigenfunction with constant sign in~$\Omega$,
see~\cite[Theorem~1.2.5]{Hen}. 
By the Faber--Krahn inequality, 
\begin{equation} \label{e:faber-krahn}
\gl_1(\Omega) \,\ge\, \gl_1(B^n(\Omega)) 
\end{equation}
where $B^n(\Omega)$ is the round ball in~$\RR^n$ with the same volume 
as~$\Omega$.
Moreover, equality holds in~\eqref{e:faber-krahn} if and only if 
$\Omega = B^n(\Omega)$, see~\cite{Faber} and~\cite{Krahn}. 
In other words, round balls are minimizers for $\gl_1$ among domains 
of the same volume.
This result can also be obtained by reasoning as follows. 
Consider the functional $\Omega \to \gl_1 (\Omega)$ for all 
smooth bounded domains $\Omega$ in $\RR^n$ of the same volume, 
say $\textnormal{Vol}(\Omega) = \ga$. 
A classical result due to Garabedian and Schiffer asserts 
that $\Omega$ is a critical point for $\gl_1$ 
(among domains of volume $\ga$) if and only if 
the first eigenfunction of the Laplacian in $\Omega$ 
with $0$~Dirichlet boundary condition has also constant Neumann data 
at the boundary, see~\cite{Gar-Schif}. 
In this case, we say that $\Omega$ is an extremal domain 
for the first eigenvalue of the Laplacian, or simply an {\it extremal domain}. 
Extremal domains are then characterized as the domains for which the 
{\it over-determined}\/ system
\begin{equation} \label{e:sys2}
\left\{\begin{array} {ll}
\Delta\, u + \gl\, u = 0 &\mbox{in }\; \Omega\\
u=0 & \mbox{on }\; \pp \Omega \\
\langle \nabla u, \nu \rangle = \const &\mbox{on }\; \pp \Omega 
\end{array}\right.
\end{equation}
has a positive solution 
(here $\nu$ is the outward unit normal vector field along $\pp \Omega$). 
By a classical result due to J.~Serrin the only domains 
for which the system~\eqref{e:sys2} has a positive solution 
are round balls, see~\cite{Serrin}. 
One then checks that round balls are minimizers.

For domains with infinite volume, at first sight 
one cannot ask for ``a domain that minimizes $\gl_1$''.
Indeed, with $c\, \Omega = \{ c\;\! z \mid z \in \Omega \}$ we have
$$
\gl_1( c\, \Omega) \,=\, c^{-2} \gl_1(\Omega) ,\quad c>0 .
$$
On the other hand, system~\eqref{e:sys2} can be studied also for 
unbounded domains. 
Therefore, it is natural to determine all domains~$\Omega$ 
for which~\eqref{e:sys2} has a positive solution. 
This is an open problem. 
We will continue to call such a domain an \textit{extremal domain}. 
In the non-compact case, this definition does not have a geometric meaning, 
except for domains which along each coordinate direction of~$\RR^n$
are bounded or periodic.
In the case of periodic directions, one obtains extremal domains for 
the first eigenvalue of the Laplacian in flat tori,
cf.~Remark~\ref{rem:torus} below. 

Berestycki, Caffarelli and Nirenberg conjectured in~\cite{BCN}
that if $f$~is a Lipschitz function on a domain $\Omega$ in~$\RR^n$ 
such that $\RR^n \backslash \overline{\Omega}$ is connected, 
then the existence of a bounded
positive solution to the more general system
\begin{equation} \label{e:sys3}
\left\{\begin{array} {ll}
\Delta u + f(u) = 0 & \mbox{in }\; \Omega\\
               u= 0 & \mbox{on }\; \pp \Omega \\
\langle \nabla u, \nu \rangle = \const &\mbox{on }\; \pp \Omega 
\end{array}\right.
\end{equation}
implies that $\Omega$ is a ball, 
or a half-space, 
or the complement of a ball,
or a generalized cylinder $B^k \times \RR^{n-k}$ 
where $B^k$ is a round ball in $\RR^k$.
In \cite{Sicbaldi}, the second author constructed a counter-example 
to this conjecture 
by showing that the cylinder $B^{n} \times \RR \subset \RR^{n+1}$ 
(for which it is easy to find a bounded positive solution to~\eqref{e:sys2})
can be perturbed to an unbounded domain whose boundary is a periodic 
hypersurface of revolution with respect to the $\RR$-axis 
and such that~\eqref{e:sys2} has a bounded positive solution. 
More precisely, 
{\it 
for each $n \ge 2$ there exists a positive number $T_* = T_*(n)$, 
a sequence of positive numbers $T_j \to T_*$,
and a sequence of non-constant $T_j$-periodic functions 
$v_j \in C^{2,\ga}(\RR)$ 
of mean~zero (over the period) that converges to~$0$ in $C^{2,\ga}(\RR)$ 
such that the domains
$$
\Omega_j \,=\, \left\{ (x,t) \in \RR^{n} \times \RR \, \mid \, \|x\| < 1+v_j(t) \right\}
$$
have a positive solution $u_j \in C^{2,\ga}(\Omega_j)$ 
to the problem~\eqref{e:sys2}.
The solution $u_j$ is $T_j$-periodic in~$t$ and hence bounded.
}

\subsection{Main results}
The goal of this paper is to show that these domains $\Omega_j$ 
(introduced in \cite{Sicbaldi} by the second author) 
belong to a smooth bifurcating family of domains, 
to determine their approximate shape for small bifurcation values,
and to determine the bifurcation values~$T_*(n)$. 
Our main result is the following.

\begin{theorem} \label{t:main}
Let $C^{2,\ga}_{\textnormal{even},0}(\RR / 2 \pi \ZZ)$ 
be the space of even $2 \pi$-periodic $C^{2,\ga}$ functions of mean zero. 
For each $n \ge 1$ there exists a positive number $T_*=T_*(n)$
and a smooth map
$$
\begin{array}{ccc}
(-\eps, \eps) & \to & C^{2,\ga}_{\textnormal{even},0}(\RR / 2 \pi \ZZ) \times \RR\\
s & \mapsto & (w_s,T_s) 
\end{array}
$$
with $w_0=0$, $T_0=T_*$ 
and such that for each $s \in (-\eps, \eps)$ the system~\eqref{e:sys2}
has a positive solution $u_s \in C^{2,\ga}(\Omega_s)$ on the modified cylinder
\begin{equation} \label{omega_s}
\Omega_s \,=\,
\left\{ (x,t) \in \RR^{n} \times \RR  \,\, \mid \,\, \|x\| < 
1+s \cos \left( \frac{2\pi}{T_s}\,t \right) + s\2 w_s \left( \frac{2\pi}{T_s}\,t \right) \right\} .
\end{equation}
The solution $u_s$ is $T_s$-periodic in~$t$ and hence bounded.
\end{theorem}

\begin{figure}[ht]
 \begin{center}
  \psfrag{t}{$t \in \RR$}
  \psfrag{x}{$x \in \RR^{n}$}
 \leavevmode\epsfbox{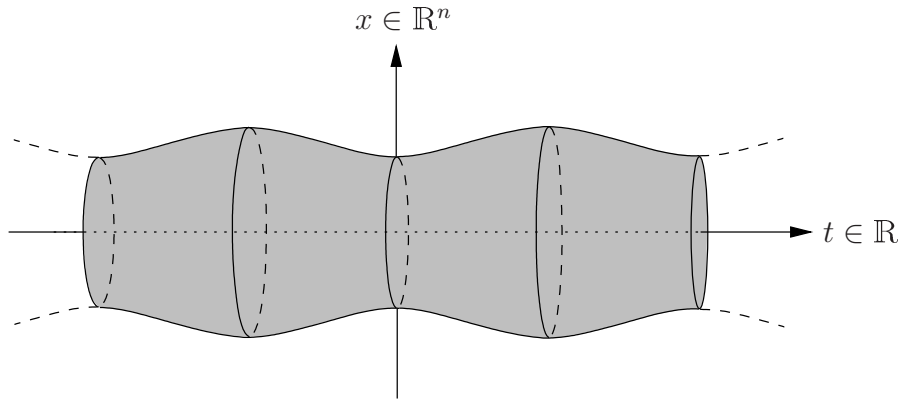}
 \end{center}
 \caption{A domain $\Omega_s$.}
 \label{figure.cyl}
\end{figure}

For $n=2$ and for $|s|$ small enough, the bifurcating domains $\Omega_s$ 
look as in Figure~\ref{figure.cyl}. 
For a figure for $n=1$ see Section~\ref{s:dim2}.

Notice that for $n=1$, the domains $\Omega_s$ do not provide counter-examples to the conjecture of 
Berestycki, Caffarelli and Nirenberg,
because $\RR^2 \setminus \Omega_s$ is not connected.

\begin{remark} \label{rem:times.r}
{\rm
From the extremal domains $\Omega_s \subset \RR^{n+1}$ and the solutions $u_s$
from Theorem~\ref{t:main} we obtain other extremal domains by adding an $\RR^k$-factor:
For each $k \ge 1$ the domains $\Omega_s^k := \Omega_s \times \RR^k$ are extremal domains
in~$\RR^{n+1+k}$ with solutions $u_s^k (x,t,y) := u_s (x,t)$ (where $y \in \RR^k$).
For instance, in~$\RR^3$ we then have the ``wavy cylinder'' in Figure~\ref{figure.cyl},
and the ``wavy board'' obtained by taking the product of the wavy band in Figure~\ref{figure.cyldim2}
with~$\RR$. 
Notice that $\RR^{n+1+k} \setminus \Omega_s^k$ is connected if and only if~$n \ge 2$.
}
\end{remark}

\begin{remark} \label{rem:torus}
{\rm
The characterization of extremal domains described in Section~\ref{ss:problem} 
more generally holds for domains in Riemannian manifolds: 
Given a Riemannian manifold $(M,g)$, a domain $\Omega \subset M$ 
of given finite volume is a critical point of 
$\Omega \to \gl_1(\Omega)$, 
where $\gl_1(\Omega)$ is the first eigenvalue of 
the Laplace--Beltrami operator~$-\Delta_g$,
if and only if the over-determined system
\begin{equation} \label{e:sys4}
\left\{\begin{array} {ll}
\Delta_g\, u + \gl\, u = 0 &\mbox{in }\; \Omega\\
u=0 & \mbox{on }\; \pp \Omega \\
g(\nabla u, \nu ) = \const &\mbox{on }\; \pp \Omega 
\end{array}\right.
\end{equation}
has a positive solution 
(here $\nu$ is the outward unit normal vector to $\pp \Omega$ 
with respect to~$g$), 
see~\cite{Elsoufi} and~\cite{Pac-Sic}. 
Theorem~\ref{t:main} thus implies that the full tori 
$$
\widetilde\Omega_s \,=\,
\left\{ (x,t) \in \RR^{n} \times \RR/T_s\, \ZZ  \;\Big|\; \|x\| < 
1+s \cos \left( \frac{2\pi}{T_s}\,t \right) + s\2 v_s \left( \frac{2\pi}{T_s}\,t \right) \right\} 
$$
are extremal domains in the manifold $\RR^{n} \times \RR / T_s \ZZ$ 
with the metric induced by the Euclidean metric.
\diam
}
\end{remark}

\ni
{\bf Open problem~1.}
{\it
Are the extremal domains $\widetilde\Omega_s$ in $\RR^{n} \times \RR/T_s \ZZ$
(local) minima for the functional $\Omega \to \gl_1(\Omega)$ ?
}

\b
It follows from our proof of Theorem~\ref{t:main} and from the Implicit Function Theorem
that the family~$\Omega_s$ is unique among those smooth families 
of extremal domains bifurcating from the straight cylinder
that are rotationally symmetric with respect to~$\RR^n$ and periodic with respect to~$\RR$.
A much stronger uniqueness property should hold.
Indeed, the existence problem of extremal domains near the solid cylinder, say in $\RR^3$, 
is tightly related to the existence problem of positive constant mean curvature surfaces near the cylinder, see Sections~\ref{s:mot} and~\ref{s:proofI}.
Any positive constant mean curvature surface with two ends (that is properly embedded and complete)
must be a Delaunay surface, 
by a result of Korevaar, Kusner, and Solomon, \cite{KKS}.
We thus ask:

\m
\ni
{\bf Open problem~2.}
{\it
Assume that $\Omega$ is an unbounded extremal domain in~$\RR^{n+1}$ 
that is contained in a solid cylinder.
Is it then true that $\Omega$ belongs to the family $\Omega_s$ ?
}

\b
We also determine the bifurcation values~$T_* = T_*(n)$.
It has been proved in \cite{Sicbaldi} that $T_*(n) < \frac{2\pi}{\sqrt{n-1}}$.
In particular, $T_*(n) \to 0$ as $n \to \infty$.
We shall show in Section~\ref{s:dim2} that $T_*(1) = 4$.
Fix now $n \ge 2$ and define $\nu = \frac{n-2}2$. Write $T_\nu$ for $T_*(n)$. 

\begin{theorem} \label{p:T}
Let $J_\nu \colon (0,+\infty) \to \RR$ be the Bessel function of the first kind.
Let $j_\nu$ be its smallest positive zero. 
Then the function $s J_{\nu-1}(s)+J_\nu(s)$ has a unique zero on the interval $(0,j_\nu)$, say $\rho_\nu$, and
$$
T_\nu \,=\, \frac{2\pi}{\sqrt{j_\nu^2-\rho_\nu^2}} .
$$
In particular,
$$
T_\nu \,=\, \sqrt{2}\2 \pi \2 \nu^{-1/2} +O(\nu^{-7/6}) .
$$
Furthermore, the sequence $T_\nu$ is strictly decreasing to~$0$.
\end{theorem}

The numbers $T_\nu$ for $\nu \le 10$ are given in Section~\ref{s:pT}. 
In particular, for $n=2, 3$ and $4$ (corresponding to the bifurcation of the straight cylinder in $\RR^3$, $\RR^4$ and $\RR^5$) the values of $T_\nu$ are
$$
T_0 \approx 3.06362, \quad T_{\frac 12} \approx 2.61931, \quad T_1 \approx 2.34104 .
$$

\ni
{\bf Open problem~3.}
{\it
Is the bifurcation at $T_*(n)$ sub-critical, critical, or super-critical\2? 
In other words, 
$\left. \partial_s(T_s) \right|_{s=0} < 0$, 
$\left. \partial_s(T_s) \right|_{s=0} = 0$, or 
$\left. \partial_s(T_s) \right|_{s=0} > 0$\2?
}

\b
The paper is organized as follows. 
In Section~\ref{s:mot} we show how the existence of Delaunay surfaces 
(i.e., constant mean curvature surfaces of revolution in~$\RR^3$ that are different from the cylinder) 
can be proved by means of a bifurcation theorem due to Crandall and Rabinowitz. 
We will follow the same line of arguments to prove Theorem~\ref{t:main}
in Sections~\ref{s:proofI} to~\ref{s:dim2}.
%
In Section~\ref{s:pT} we prove Theorem~\ref{p:T} on the bifurcation 
values~$T_*(n)$. 

\b
\ni
{\bf Acknowledgments.}
Most of this paper was written in June~2010, when the second author visited
Universit\'e de Neuch\^atel. 
The second author is grateful to Bruno Colbois and Alexandre Girouard for
their warm hospitality.
The first author thanks FRUMAM for its hospitality during the workshop 
``Probl\`emes aux valeurs propres et probl\`emes surd\'etermin\'es'' 
at Marseille in December~2010.
We both thank Frank Pacard for helpful discussions and for 
kindly allowing us to include the exposition in Section~\ref{s:mot}.


\section{The Delaunay surface via
the Crandall--Rabinowitz Theorem}  \label{s:mot}

\ni
Our proof of Theorem~\ref{t:main} is motivated by the following argument that proves 
the existence of Delaunay surfaces by means of the Crandall--Rabinowitz bifurcation theorem.
The material of this section was explained by Frank Pacard to the second author
when he was his PhD student.

\m
We start with some generalities. 
Let $\Sigma$ be an embedded hypersurface in $\RR^{n+1}$ of codimension~1. 
We denote by~$II$ its second fundamental form defined by
$$
II(X, Y) \,=\, -\left\langle \nabla_X\, N, Y\right\rangle
$$
for all vector fields $X,Y$ in the tangent bundle $T\,\Sigma$. 
Here~$N$ is the unit normal vector field on~$\Sigma$, and 
$\left\langle \cdot, \cdot \right\rangle$ denotes the standard 
scalar product of $\RR^{n+1}$. 
The mean curvature~$H$ of~$\Sigma$ is defined to be the average 
of the principal curvatures, i.e.~of the eigenvalues $k_1, \dots , k_{n}$ 
of the shape operator $A \colon T\, \Sigma \longrightarrow T\,\Sigma$ 
given by the endomorphism
$$
\left\langle A\,X, Y\right\rangle \,=\, -II(X, Y) .
$$
Hence
$$
H(\Sigma) \,=\, \frac{1}{n}\, \sum_{i=1}^n k_i .
$$
Given a sufficiently smooth function $w$ defined on $\Sigma$ we can define 
the normal graph $\Sigma_w$ of~$w$ over~$\Sigma$,
$$
\Sigma_w \,=\, 
\left\{ p + w(p)\, N(p) \in \RR^{n+1} \2 \mid \, p \in \Sigma \right\} ,
$$
and consider the operator $w \mapsto H(\Sigma_w)$ that associates to~$w$
the mean curvature of $\Sigma_w$. 
The linearization of this operator at $w=0$ is given by the Jacobi operator:
$$
D_{w} \left.H(\Sigma_w)\right|_{w=0} \,=\, 
\frac{1}{n}\, \left(\Delta_g + \sum_{i=1}^n k_i^2\right),
$$ 
where $g$ the metric induced on $\Sigma$ by the Euclidean metric and 
$-\Delta_g$ is the Laplace--Beltrami operator on~$\Sigma$. 
All these facts are well-known, 
and we refer to~\cite{cold-min} for further details.

\m
In 1841, C.~Delaunay discovered a beautiful one-parameter family 
of complete, embedded, non-compact surfaces $D_\gs$ in~$\RR^3$, 
$\gs>0$, whose mean curvature is constant, see \cite{delaunay}. 
These surfaces are invariant under rotation about an axis and periodic in 
the direction of this axis. 
The Delaunay surface~$D_\gs$ can be parametrized by 
$$
X_\gs(\theta,t) \,=\, \bigl( y(t) \cos \theta, y(t) \sin \theta, z(t) \bigr)
$$
for $(\theta,t) \in S^1 \times \RR$, where the function~$y$ is the smooth solution of
$$
(y'(t))^2 \,=\, y^2(t) - \left( \frac{y^2(t) + \gs}{2}\right)^2
$$
and $z$ is the solution (up to a constant) of 
$$
z'(t) \,=\, \left( \frac{y^2(t) + \gs}{2}\right).
$$
When $\gs =1$, the Delaunay surface is nothing but the cylinder 
$D_1 = S^1 \times \RR$. 
It is easy to compute the mean curvature of the family~$D_\gs$ 
and to check that it is equal to~1 for all~$\gs$. 
One can obtain each Delaunay surface~$D_\gs$ 
by taking the surface of revolution generated by the roulette of an ellipse, 
i.e.~the trace of a focus of an ellipse~$\ell$ as $\ell$ rolls 
along a straight line in the plane. 
In particular, these surfaces are periodic in the direction of 
the axis of revolution. 
When the ellipse~$\ell$ degenerates to a circle, 
the roulette of~$\ell$ becomes a straight line and generates 
the straight cylinder, and when $\gs \to 0$, $D_\gs$ tends to the singular surface
which is the union of infinitely many spheres of radius~$1/2$ centred at the points
$(0,0,n)$, $n \in \ZZ$.
For further details about this geometric description of Delaunay surfaces we refer to~\cite{eells}.

\m
We now prove the existence of Delaunay surfaces by a bifurcation argument,
using a bifurcation theorem due to M.~Crandall and P.~Rabinowitz. 
Their theorem applies to Delaunay surfaces in a simple way.
We shall use the same method to prove Theorem~\ref{t:main}. 
The phenomenon underlying our existence proof of Delaunay surfaces
is the Plateau--Rayleigh instability of the cylinder, \cite{Ray}.

Consider the straight cylinder of radius~1, in cylindrical coordinates:
$$
C_{1} \,=\, \{(\rho, \theta,t) \in (0,+\infty) \times S^1 \times \RR \,\, \mid \,\, \rho = 1 \} .
$$
Let $w$ be a $C^2$-function~on $S^1 \times \RR/2 \pi \ZZ$.
In Fourier series,
$$
w (\theta,t) \,=\, 
\sum_{j,k \ge 0} 
\bigl(\ga_j\,\cos(j\, \theta) + \gb_j\,\sin(j\,\theta) \bigr) \, 
\bigl(a_k\,\cos\left(k\, t\right) + b_k\, \sin \left(k\, t\right) \bigr) .
$$
%
If $w(\theta,t) > -1$ for all $\theta,t$, we consider, 
for each $T>0$, 
the normal graph $C_{1+w}^T$ over the cylinder~$C_1$ of~$w$ rescaled to period~$T$, 
$$
C_{1+w}^T \,:=\,  
\left\{
(\rho, \theta,t) \in (0,+\infty) \times S^1 \times \RR \,\, \mid \,\,
\rho = 1 + w \left( \theta, \frac{2\pi}{T}\,t \right)
\right\} .
$$
Define the operator
$$
\widetilde F(w,T) \,=\, 1 - H \left(C_{1+w}^T \right)
$$
where $H$ is the mean curvature. 
Then $\widetilde F(w,T)$ is a function on $S^1 \times \RR$ 
of period~$T$ in the second variable. 
Therefore,
\begin{equation}\label{operator}
F(w, T)\, (\theta, t) \,:=\, 
         \widetilde F(w, T) \left(\theta, \frac{T}{2\pi}\,t\right)
\end{equation}
is a function on $S^1 \times \RR/2 \pi \ZZ$. 
Note that $F(0,T) = 0$ for all $T>0$, 
because for $w=0$ the surface $C_{1+w}^T$ is the cylinder~$C_1$ whose mean curvature is~1. 
If we found a non-trivial solution $(w,T)$ of the equation $F(w,T)=0$, 
we would obtain a 
constant mean curvature surface different from~$C_1$.
%
In order to solve this equation, we consider the linearization 
of the operator~$F$ with respect to~$w$ and computed at $(w,T) = (0,T)$. 
As mentioned above, the linearization of the mean curvature operator 
for normal graphs over a given surface with respect to~$w$ 
computed at~$w=0$ is the Jacobi operator.
Since the Laplace--Beltrami operator on~$C_1$ 
(with the metric induced by the Euclidean metric)
is $-\partial^2_{\theta} - \partial^2_t$, and since 
the principal curvatures~$k_i$ of $C_1$ are equal to~$0$ and~$1$, we find that 
$$
D_w\,F(0,T) \,=\, 
- \frac{1}{2}\, \left(\partial^2_{\theta} + \left( \frac{2\pi}{T}\right)^2 \partial^2_t +1\right) .
$$
For each $j,k \in \NN \cup \{0\}$ and each $T>0$, 
the four $1$-dimensional spaces generated by the functions
$$
\cos (j \2 \theta) \2 \cos (k \2 t) , \quad 
\cos (j \2 \theta) \2 \sin (k \2 t) , \quad
\sin (j \2 \theta) \2 \cos (k \2 t) , \quad 
\sin (j \2 \theta) \2 \sin (k \2 t)
$$ 
are eigenspaces of $D_w\,F(0,T)$ with eigenvalue 
$$
\sigma_{j,k}(T) \,=\, \frac{1}{2}\, \left(j^2 - 1 + \left(\frac{2 \pi k}{T}\right)^2\right) .
$$
Clearly, 
\begin{itemize}
  \item $\sigma_{j,k}(T) \neq 0$ for all $T>0$ \, if $j \ge 2$, or if $j =1$ and $k\ge 1$;
	\item $\sigma_{1,0}(T) =0$ for all $T>0$;
	\item $\sigma_{0,k}(T) = 0$ only for $T = 2 \pi k$ and $k \ge 1$; 
	      moreover $\sigma_{0,k}(T)$ changes sign at these points.
\end{itemize}
It follows that $\Ker D_wF (0,T)$ is $2$-dimensional (spanned by $\cos \theta$, $\sin \theta$)
if $T>0$ and $T \notin 2\pi \NN$, 
and that $\Ker D_wF (0,T)$ is $4$-dimensional 
(spanned by $\cos \theta$, $\sin \theta$, $\cos (kt)$, $\sin (kt)$) if $T \in 2\pi \NN$.

We will now bring into play an abstract bifurcation theorem,
which is due to Crandall and Rabinowitz. 
For the proof and for many other applications we refer to~\cite{Kielhofer, Smoller} 
and to the original exposition~\cite{crand-rab}.

\begin{theorem}  \label{t:CR}
{\rm \bf (Crandall--Rabinowitz Bifurcation Theorem)}
Let $X$ and $Y$ be Banach spaces, and let $U \subset X$ and $\Lambda \subset \RR$
be open subsets, where we assume $0 \in U$. 
Denote the elements of $U$ by~$w$ and the elements of $\Lambda$ by~$T$.
Let $F \colon U \times \Lambda \to Y$ be a $C^\infty$-smooth function such that

\s
\begin{itemize}
\item[i)]
$F(0,T) =0$ for all $T \in \Lambda$;

\s
\item[ii)]
$\Ker D_w\, F(0,T_0) = \RR \2 w_0$\,  for some $T_0 \in \Lambda$ and some $w_0 \in X \setminus \{0\}$;

\s
\item[iii)]
$\codim \Im D_w\,F(0,T_0) =1$;

\s
\item[iv)]
$D_T D_w\, F(0,T_0) (w_0) \notin \Im D_w\, F(0,T_0)$.
\end{itemize}
Choose a linear subspace $\dot X \subset X$ such that $\RR \2 w_0 \oplus \dot X = X$.
Then there exists a $C^\infty$-smooth curve 
$$
(-\eps, \eps) \to \dot X \times \RR, 
\quad
s \mapsto \left( w(s), T(s) \right) 
$$
such that 

\s
\begin{itemize}
\item[1)]
$w(0) = 0$ and $T(0) = T_0$;

\s
\item[2)]
$s \left( w_0+w(s) \right) \in U$ and $T(s) \in \RR$;

\s
\item[3)]
$F \bigl( s \left( w_0+w(s) \right), T(s) \bigr) =0$.
\end{itemize}

\s \ni
Moreover, there is a neighbourhood $\cn$ of $(0,T_0) \in X \times \RR$ 
such that $\left\{ s \left( w_0+w(s) \right), T(s) \bigr) \right\}$ is the 
only branch in~$\cn$ that bifurcates from $\{(0,T) \mid T \in \Lambda \}$.
\end{theorem}

The theorem is useful for finding non-trivial solution of an equation 
$F(x,\lambda)=0$, where $x$ belongs to a Banach space  
and $\lambda$ is a real number. 
It says that under the given hypothesis, there is a smooth bifurcation 
into the direction of the kernel of $D_wF$ for the solution of 
$F(x, \lambda)=0$, and that there is no other nearby bifurcation.

\s 
In order to apply Theorem~\ref{t:CR},
we now restrict the operator~$F$ defined in~\eqref{operator} to functions 
that are independent of~$\theta$
(so as to get rid of the functions $\cos \theta$, $\sin \theta$ in the kernel of
$D_w F(0,T)$) 
and that are even (so as to have a 1-dimensional kernel for $T \in 2\pi \NN$).
We can also assume that the functions $w$ have zero mean.
In other words, we look for new constant mean curvature surfaces among deformations of~$C_1$
that are surfaces of revolution, even in the $t$-direction.
We hence consider the Banach space 
$$
X \,=\, C^{2,\alpha}_{\textnormal{even},0} (\RR/2 \pi \ZZ)
$$
of even $2\pi$-periodic functions of zero mean whose second derivative is H\"older continuous.
Moreover, define the open subset $U = \left\{ w \in X \mid w(t) >-1 \mbox{ for all } t \right\}$ 
of~$X$,
and the Banach space
$$
Y \,=\, C^{0,\alpha}_{\textnormal{even},0} (\RR/2 \pi \ZZ) .
$$
Furthermore, chose $\Lambda = (0,+\infty) \subset \RR$.
Then the operator $F$ defined as above restricts to the operator
$$
F \colon U \times \Lambda \,\to\, Y .
$$
With
$$
\sigma_k(T) \,:=\, \sigma_{0,k}(T) \,=\, 
\frac{1}{2}\, \left(- 1 + \left(\frac{2 \pi k}{T}\right)^2\right) ,
$$
its linearization with respect to $w$ at $T_0 := 2\pi$ is
$$
D_w\, F(0,T_0) \left( \sum_{k \ge 1} a_k\, \cos(k\,t) \right) \,=\, 
\sum_{k \ge 1} \sigma_k(T_0)\, a_k \, \cos(k\,t) \,=\,
\sum_{k \ge 1} \frac 12 (k^2-1)\, a_k \, \cos(k\,t) .
$$
Hence,
$$
\Ker D_w\, F(0,T_0) = \RR\, \cos \1 t .
$$
Moreover, the image $\Im D_w\, F(0,T_0)$ is the closure of 
$\bigoplus_{k \geq 2} \RR \cos (kt)$ in~$Y$;
its complement in $Y$ is the $1$-dimensional space spanned by $\cos \1 t$.
Finally,
$$
D_T \1 D_w\1 F(0,T_0) \left( \cos \2 t \right) \,=\,
\left.\frac{\pp \sigma_1(T)}{\pp T} \right|_{T=T_0} \,\cos \2 t \,=\, 
- \frac{1}{2\pi} \cos \1 t  \,\notin\, \Im D_w \1 F(0,T_0).
$$
With $w_0 = \cos \1 t$ and $\dot X$ the closure of $\bigoplus_{k \geq 2} \RR \cos (kt)$ in~$X$, 
the Crandall--Rabinowitz bifurcation theorem applies and yields the existence of $C^\infty$-smooth curve 
$$
(-\eps, \eps) \,\to\, \dot X \times \RR, 
\quad
s \,\mapsto\, \left( w(s), T(s) \right) 
$$
such that 
\begin{itemize}
\item[1)]
$w(0) = 0$ and $T(0) = T_0$;

\s
\item[2)]
$F \bigl( s \left( w_0+w(s) \right), T(s) \bigr) =0$,
\end{itemize}
i.e.~(by the definition of the operator~$F$) the existence of a $C^\infty$-smooth family of surfaces of revolution that have mean curvature constant and equal to 1, 
bifurcating from the cylinder~$C_1$. 
That these surfaces are Delaunay surfaces 
follows from Sturm's variational characterization of constant mean curvature surfaces 
of revolution,~\cite{delaunay,eells}. 

\begin{remark}
{\rm
The boundaries of the new domains~$\Omega_s \subset \RR^3$ 
described in Theorem~\ref{t:main} are not Delaunay surfaces
(at least not for $|s|$ small).
Indeed, Delaunay surfaces bifurcate from the cylinder at $T_0=2\pi$,
while the domains $\Omega_s$ bifurcate from the cylinder at $T_*(2) \approx 3.06362$.  
}
\end{remark}


\section{Rephrasing the problem for extremal domains}  \label{s:proofI}

\ni
We want to follow the proof of the existence of Delaunay surfaces given in the previous section 
in order to prove the existence of a smooth family of normal graphs over the straight cylinder 
such that the first eigenfunction of the Dirichlet Laplacian has constant Neumann data. 
In this section we recall the set-up from~\cite{Sicbaldi},
where the second author studied the Dirichlet-to-Neumann
operator
that associates to a periodic function~$v$
the normal derivative of the first eigenfunction of the domain defined by the  
normal graph of~$v$ over the straight cylinder, 
and computed the linearization of this operator. 
The novelty of this paper is the analysis of the kernel of the linearized operator;
it will be carried out in Sections~\ref{s:Bessel} to \ref{s:extremalviaCR}.

\s
The manifold $\RR/ 2\pi \ZZ$ will always be considered with the metric 
induced by the Euclidean metric. 
Motivated by the previous section, we consider the Banach space
$\mathcal C^{2,\alpha}_{\textnormal{even},0} (\mathbb{R}/2\pi \mathbb{Z})$ 
of even functions on $\mathbb{R}/2\pi \mathbb{Z}$ of mean~$0$. 
For each function $v \in \mathcal C^{2,\alpha}_{\textnormal{even},0} (\mathbb{R}/2\pi \mathbb{Z})$ 
with $v(t) >-1$ for all $t$, the domain
$$
C_{1+v}^T  \,:=\,  
\left\{  
(x,t) \in \RR^n \times \RR/T \mathbb{Z}  \,\, \mid  \,\, 
0 \le \|x\| < 1+ v \left(\frac{2\,\pi}{T}\,t \right)
\right\} 
$$
is well-defined for all $T>0$. 
The domain $C_{1+v}^T$ is relatively compact.
According to standard results on the Dirichlet eigenvalue problem (see \cite{gilbarg}),
there exist, for each $T>0$, a unique positive function 
$$
\phi = \phi_{v,T} \,\in\, {\mathcal C}^{2, \alpha} \left( C_{1+v}^T \right)
$$
and a constant $\lambda = \lambda_{v,T} \in \RR$ 
such that $\phi$ is a solution to the problem
\begin{equation}\label{formula}
\left\{
\begin{array}{rcccl}
	\Delta \, \phi + \lambda \,\phi &=& 0 & \textnormal{in} & C_{1+v}^T \\[1mm]
	\phi & = & 0 & \textnormal{on} & \partial C_{1+v}^T
\end{array}
\right.
\end{equation}
which is normalized by 
\begin{equation}
\int_{C_{1+v}^{2\pi}} \left(\phi \left(x,\frac{T}{2\pi}\,t \right)\right)^2\, {\rm dvol} \,=\, 1 .
\label{noral}
\end{equation}
Furthermore, $\phi$ and $\lambda$ depend smoothly on~$v$. 
We denote $\phi_1 := \phi_{0,T}$ and $\lambda_1:=\lambda_{0,T}$. 
Notice that $\phi_1$ does not depend on the $t$~variable and is radial in the $x$~variable. 
(Indeed, $\phi_1$ is nothing but the first eigenfunction of the Dirichlet Laplacian 
over the unit ball~$\mathbb{B}^n$ in~$\RR^n$ normalized to have $L^2$-norm $\frac{1}{2\pi}$.)
We can thus consider $\phi_1$ as a function of $r:=\|x\|$, and we write 
\begin{equation} \label{1:eigenf}
\varphi_1(r) = \phi_1(x).
\end{equation}
We define the Dirichlet-to-Neumann
operator
$$
\widetilde F(v,T)  \,=\,  
\langle \nabla \phi , \nu \rangle  \, |_{\pp C_{1+v}^{T}}  
- \frac{1}{{\rm Vol} (\pp C_{1+v}^{T})} \, \int_{\pp C_{1+v}^{T}} \langle \nabla \phi , \nu \rangle \, \mbox{dvol} \, ,
$$
where $\nu$ denotes the unit normal vector field on $\pp C_{1+v}^{T}$ 
and where~$\phi = \phi_{v,T}$ is the solution of~\eqref{formula}. 
The function 
$$
\widetilde F (v,T) \colon \pp C_{1+v}^T \cong \partial (\mathbb{B}^{n}) \times \RR/T \ZZ \,\to\, \RR 
$$ 
depends only on the variable $t \in \RR/T \ZZ$,
since $v$ has this property.
It is an even function; indeed, $v$ is even, and hence $\phi_{v,T}$ is
even, since the first eigenvalue $\lambda_{v,T}$ is simple.  
Moreover, $\widetilde F(v,T)$ has mean~$0$. 
We rescale $\widetilde F$ and define
$$
F(v,T)\, (t) \,=\, \widetilde F (v,T)\, \left( \frac{T}{2\pi}\,t \right) .
$$
Schauder's estimates imply that~$F$ takes values in 
$\mathcal C^{1,\alpha}_{\textnormal{even},0} (\RR/2\pi \ZZ)$.
With 
$$
U \,:=\, 
\left\{ v \in \mathcal C^{2,\alpha}_{\textnormal{even},0} (\RR/2\pi \ZZ) \mid v(t) >-1 \mbox{ for all } t \right\}
$$ 
we thus have
$$
F \colon U \times (0,+\infty) \,\to\, \mathcal C^{1,\alpha}_{\textnormal{even},0} (\RR/2\pi \ZZ) .
$$
Also notice that $F(0,T)=0$ for all $T>0$, and that $F$ is smooth.

The following result is proved in~\cite{Sicbaldi}.
\begin{proposition} \label{linearization}
The linearized operator 
$$
H_T := D_w \2 F(0,T) \colon 
\mathcal C^{2,\alpha}_{\textnormal{even},0} (\RR/2\pi \ZZ) \longrightarrow 
\mathcal C^{1,\alpha}_{\textnormal{even},0} (\RR/2\pi \ZZ)  
$$
is a formally self adjoint, first order elliptic operator. It preserves the eigenspaces 
\[
V_k = \RR \, \cos (kt)  
\]
for all~$k$ and all~$T>0$, and we have
\begin{equation}\label{def:op:H}
H_T(w)(t) \,=\, 
\left.
\left(\partial_r \psi + \partial^2_r \phi_1 \cdot w \left( \frac{2\pi t}{T}\right) \right)\right|_{\pp C_1^T}
\end{equation}
where $\psi$ is the unique solution of 
\begin{equation} \label{def:psi}
\left\{\begin{array} {ll}
\Delta \psi + \lambda_1 \psi = 0 & \mbox{in }\; C_1^T\\
               \psi = -\partial_r \phi_1 \cdot w(2\pi t/T) & \mbox{on }\; \pp C_1^T
\end{array}\right.
\end{equation}
which is $L^2(C_1^T)$-orthogonal to $\phi_1$, and where $r = \|x\|$.
\end{proposition}
Write 
$$
w(t) \,=\, \sum_{k \ge 1} a_k\2 \cos(k\1t).
$$
Since $H_T$ preserves the eigenspaces, 
$$
H_T(w)(t) \,=\, \sum_{k \ge 1} \sigma_{k}(T) \2 a_k \2 \cos(k\1 t).
$$
We use \eqref{def:op:H} and \eqref{def:psi} to describe $\sigma_k(T)$ 
as the solution of an ordinary differential equation: 
The solution~$\psi$ of~\eqref{def:psi} is differentiable, 
and even with respect to~$x$ for fixed~$t$. 
Therefore, for each~$t$, the derivative of~$\psi$ with respect to~$r$ vanishes at 0: 
$\partial_r \psi|_{r=0} = 0$. 
Hence,
\begin{equation} \label{sigma:1:def}
\sigma_k(T) \,=\, c_k'(1) + \varphi_1''(1)
\end{equation}
where for $n \geq 2$, $c_k$ is the continuous solution on $[0,1]$ 
of the ordinary differential equation
$$
\left( \pp_r^2 + \frac{n-1}{r}\, \pp_r + \lambda_1 - \left(\frac{2\,\pi\,k}{T}\right)^2\right)\, c_k \,=\, 0
$$
such that $c_k(1) = -\varphi_1'(1)$, 
while for $n=1$, $c_k$ is the solution on $[0,1]$ of the ordinary differential equation
$$
\left( \pp_r^2 + \lambda_1 - \left(\frac{2\,\pi\,k}{T}\right)^2\right)\, c_k \,=\, 0
$$
such that $c_k(1) = -\varphi_1'(1)$ and $c_k'(0)=0$.
Notice that for all $k \geq 1$ and all $n \geq 1$
$$
\sigma_k(T) \,=\, \sigma_1 \left( \frac{T}{k} \right) .
$$ 
Our next aim is to find an explicit expression for the function $\sigma_1$ 
in order to describe the spectrum of the linearized operator, 
to read off its kernel, and to find the codimension of its image.
We first consider the case $n \geq 2$, for which we need Bessel functions.
The case $n=1$ is discussed in Section~\ref{s:dim2}.


\section{Recollection on Bessel functions}  \label{s:Bessel}

\ni
In what follows we shall use several basic properties of Bessel functions.  
For the readers convenience, we recall the definition of the Bessel functions 
$J_\tau$ and $I_\tau$, and state their principal properties. 
For proofs we refer to~\cite[Ch.~III]{Watson}.

\subsection{The functions $J_\tau$}
For $\tau \ge 0$ the Bessel function of the first kind 
$J_\tau \colon \RR \to \RR$ is the solution of 
the differential equation
\begin{equation} \label{diff:eq:J}
s^2\, y''(s) + s\,y'(s) + (s^2 - \tau^2)\,y(s) \,=\, 0
\end{equation}
whose power series expansion is 
\begin{equation}  \label{exp:J}
J_\tau (s) \,=\, \sum_{m=0}^\infty \frac{(-1)^m(\frac 12 s)^{\tau+2m}}{m!\,\Gamma (\tau+m+1)} .
\end{equation}
We read off that 
\begin{equation} \label{eJ:at0}
J_0(0)=1, \quad J_\tau(0) = 0 \; \mbox{ for all } \1 \tau > 0 . 
\end{equation}
The power series~\eqref{exp:J} defines a solution $J_\tau \colon (0,\infty) \to \RR$ of~\eqref{diff:eq:J}
also for $\tau < 0$.
If $\tau = n$ is an integer, then 
\[
J_{-n}(s) \,=\, (-1)^n\, J_n(s)
\]
and $J_n$ is bounded near~$0$.
If $\tau$ is not an integer, then the function $J_\tau(s)$ is bounded near~$0$ if $\tau >0$ 
but diverges as $s \to 0$ if $\tau < 0$. 
The functions $J_\tau(s)$ and $J_{-\tau}(s)$ are therefore linearly independent, 
and hence are the two solutions of the differential equation~\eqref{diff:eq:J} on~$(0,\infty)$.

For all $\tau \in \RR$ and all $s>0$ we have the recurrence relations
\begin{eqnarray}
J_{\tau-1}(s)+J_{\tau+1}(s) &=& \frac{2\tau}{s} J_\tau(s) , \label{e:J1} \\
J_{\tau-1}(s)-J_{\tau+1}(s) &=& 2 J_\tau'(s) , \label{e:J2} \\
sJ_{\tau}'(s)+\tau J_{\tau}(s) &=& \phantom{-}s J_{\tau-1}(s) , \label{e:J3} \\
sJ_{\tau}'(s)-\tau J_{\tau}(s) &=& -s J_{\tau+1}(s) . \label{e:J4} 
\end{eqnarray}
Another important property that we will use often is that the first eigenvalue~$\lambda_1$ 
of the Dirichlet Laplacian on the unit ball of~$\RR^n$, $n \ge 2$, 
is equal to the square of the first positive zero of~$J_\nu$ for $\nu = \frac{n-2}{2}$. 
Notice that $\lambda_1$ depends on~$n$. 
Moreover, the function $J_\nu$ is positive on the interval $(0,\sqrt{\lambda_1})$, 
and $J'_\nu (\sqrt{\lambda_1}) < 0$.

\subsection{The functions $I_\tau$}
For $\tau \in \RR$ the modified Bessel function of the first kind 
$I_\tau \colon \RR \to \RR$ is the solution of the differential equation
$$
s^2\, y''(s) + s\,y'(s) - (s^2+\tau^2)\,y(s) \,=\, 0
$$
whose power series expansion is
\begin{equation} \label{e:power.I}
I_\tau (s) \,=\, \sum_{m=0}^\infty \frac{(\frac 12 s)^{\tau+2m}}{m!\,\Gamma (\tau+m+1)} .
\end{equation}
We read off that $I_\tau(s) >0$ for all $\tau \in \RR$ and $s>0$, and that
\begin{equation} \label{e:at0}
I_0(0)=1, \quad I_\tau(0) = 0 \; \mbox{ for all } \1 \tau > 0. 
\end{equation}
Comparing coefficients readily shows that 
for all $\tau \in \RR$ and all $s>0$ we have the recurrence relations
\begin{eqnarray}
I_{\tau-1}(s)-I_{\tau+1}(s) &=& \frac{2\tau}{s} I_\tau(s) , \label{e:I1} \\
I_{\tau-1}(s)+I_{\tau+1}(s) &=& 2 I_\tau'(s) , \label{e:I2} \\
sI_{\tau}'(s)+\tau I_{\tau}(s) &=& s I_{\tau-1}(s) , \label{e:I3} \\
sI_{\tau}'(s)-\tau I_{\tau}(s) &=& s I_{\tau+1}(s) . \label{e:I4} 
\end{eqnarray}
We shall also make use of the asymptotics
\begin{equation} \label{asym:I}
\lim_{s \to \infty} \frac{I_\tau(s)}{\frac 1{\sqrt{2\pi s}} \1 e^s} \,=\, 1.
\end{equation}


\section{A formula for $\gs_1(T)$ when $n \ge 2$}   \label{s:proofII}

\ni
In this section we begin our analysis of the first eigenvalue $\sigma_1(T)$ of 
the linearized operator~$H_T$. 
We assume that $n \ge 2$ throughout. 
To simplify the notation, we denote the previously defined function~$c_1$ by~$c$.
Recall that for $n \geq 2$, 
$$
\sigma_1(T) \,=\, c'(1) + \varphi_1''(1)
$$
where $c$ is the continuous solution on $[0,1]$ of the ordinary differential equation
\begin{equation} \label{e:diffc}
\left(
\partial_r^2 + \frac{n-1}{r}\, \partial_r + \lambda_1 - \left(\frac{2\,\pi}{T}\right)^2
\right)\, c \,=\, 0
\end{equation}
such that $c(1) = -\varphi_1'(1)$. 
We shall distinguish three cases, according to whether the term 
$$
\lambda_1 - \left(\frac{2\,\pi}{T}\right)^2
$$ 
is negative, zero or positive. 
Recall that $\lambda_1$ depends on~$n$. 
In order to simplify notation, we put $\nu = \frac{n-2}{2}$ 
and write $\lambda_\nu$ for $\lambda_1=\lambda_1(n)$. 
As mentioned in the previous section, $\sqrt{\gl_\nu}$ is the first zero 
of~$J_\nu$. Denote
$$
j_\nu = \sqrt{\gl_\nu} 
$$
and $\mu = \frac{2\pi}{j_\nu}$. 
We shall find an explicit expression for $\sigma_1(T)$.
For $T>0$ denote 
\begin{equation} \label{sigma:left:right}
\sigma_1(T) = \left\{\begin{array}{lll}
\gs_{\lef}(T)  & \textnormal{if} & T < \mu, \\
\sigma_1(\mu)  & \textnormal{if} & T = \mu, \\
\gs_{\righ}(T) & \textnormal{if} & T > \mu.
\end{array}\right.
\end{equation}

\subsection{A formula for $\gs_{\lef}$.}
Assume that $T<\mu$. This allows us to define
\begin{equation} \label{def:xi}
\xi \,=\, \sqrt{\left(\frac{2\,\pi}{T}\right)^2-\lambda_\nu}\, .
\end{equation}
We rescale the function~$c$ by defining 
$$
\tilde c(s) = c\left(\frac{s}{\xi}\right) .
$$
In view of \eqref{e:diffc}, 
$\tilde c$ is the continuous solution on $[0,\xi]$ of 
$$ 
\left( \partial_s^2+ \frac{n-1}{s}\, \partial_s - 1\right)\,\tilde c\, = \, 0
$$
with $\tilde c(\xi) = -\varphi_1'(1)$. This equation is very similar to a modified Bessel equation. 
In order to obtain exactly a modified Bessel equation, we define the function $\hat c$ by
$$
\tilde c(s) = s^{-\nu}\, \hat c(s) .
$$
Note that $-\nu \le 0$ because $n \ge 2$. 
Hence $\hat c$ is the continuous solution on $[0,\xi]$ of
$$
\left[\partial_s^2+ \frac{1}{s}\, \partial_s - \left(1 + \frac{\nu^2}{s^2}\right)\right]\,\hat c \,=\, 0
$$
with $\hat c(\xi) = -\xi^{\nu}\, \varphi_1'(1)$. 
The solution of this ordinary differential equation is given by 
$\alpha\,I_{\nu}(s)$, where the constant $\alpha$ (depending on $\nu$ and $T$) 
is chosen such that 
$$
\alpha\, I_\nu(\xi) = -\xi^{\nu}\, \varphi_1'(1).
$$
Returning to the function~$c$, we get
$$
c(r) \,=\,  -\frac{\varphi_1'(1)}{I_\nu(\xi)}\, r^{-\nu}\, I_\nu(\xi\, r)
$$
and from \eqref{sigma:1:def} and \eqref{sigma:left:right}, using the identities~\eqref{e:I2}, \eqref{e:I3} and~\eqref{e:I4}, we obtain
\begin{eqnarray}
\gs_{\lef}(T) 
&=&
-\gf_1'(1) \frac{1}{I_\nu (\xi)} \frac 12 
\Bigl( \xi I_{\nu-1}(\xi) -2\nu I_\nu(\xi) + \xi I_{\nu+1}(\xi) \Bigr)  + \gf_1''(1) \notag \\
&=&
\gf_1''(1) - \gf_1'(1) \frac{\xi I_{\nu+1}(\xi)}{I_\nu(\xi)} .  \label{e:sl}
\end{eqnarray}
To better understand $\gs_{\lef}(T)$ we shall need the values $\gf'(1)$ and $\gf''(1)$.
From~\eqref{1:eigenf} and the definition of~$\phi_1$ we have that $\varphi_1$ is the continuous solution on $[0,1]$ of 
\[ 
\left( \partial_r^2+ \frac{n-1}{r}\, \partial_r + \lambda_\nu \right)\,\varphi_1\, = \, 0
\]
such that $\varphi_1(1)=0$, with normalization
\[
\int_0^1 \varphi_1^2(r)\, \textnormal{d}r = \frac{1}{2 \pi\, \textnormal{Vol}(S^{n-1})} .
\]
We rescale the function $\varphi_1$ and define
$$
\tilde \varphi_1(s) = \varphi_1 \left(\frac{s}{j_\nu}\right) .
$$
Hence, $\tilde \varphi_1$ is the continuous solution on $[0,j_\nu]$ of 
\begin{equation} \label{e:2s0}
\left( \partial_s^2+ \frac{n-1}{s}\, \partial_s + 1\right)\,\tilde \varphi_1 \, = \, 0
\end{equation}
with $\tilde \varphi_1 (j_\nu) = 0$ and  normalization
$$
\int_0^{j_\nu} \tilde \varphi_1^2(s)\, \textnormal{d}s = \frac{j_\nu}{2 \pi\, \textnormal{Vol}(S^{n-1})} .
$$ 
Equation~\eqref{e:2s0} is very similar to a Bessel equation. 
In order to obtain exactly a Bessel equation, we define the function~$\hat \gf_1$
by
\[
\tilde \varphi_1(s) = s^{-\nu}\, \hat \varphi_1(s) .
\]
Since $-\nu \le 0$ because $n \ge 2$, we get that $\hat \varphi_1$ is the continuous solution on $[0,j_\nu]$ of
\[ 
\left[\partial_s^2+ \frac{1}{s}\, \partial_s + \left(1 - \frac{\nu^2}{s^2}\right)\right]\,\hat \varphi_1 \,=\, 0
\]
with $\hat \varphi_1 (j_\nu) = 0$ and normalization
$$
\int_0^{j_\nu} s^{2-n}\, \hat \varphi_1^2(s)\, \textnormal{d}s = \frac{j_\nu}{2 \pi\, \textnormal{Vol}(S^{n-1})} .
$$
The solution of this ordinary differential equation is 
$\kappa_n \, J_{\nu}(s)$, where the constant $\kappa_n$ is chosen such that 
$$
\int_0^{j_\nu} \kappa_n^2\, s^{2-n}\, J_\nu^2(s)\, \textnormal{d}s 
\,=\, \frac{j_\nu}{2 \pi\, \textnormal{Vol}(S^{n-1})} .
$$
Returning to the function $\varphi_1$, we get
$$
\varphi_1 (r) = \kappa_n\, j_\nu^{-\nu}\, r^{-\nu}\, J_\nu (j_\nu\, r).
$$
It follows that
$$
\varphi_1'(r) \,=\, \kappa_n\, j_\nu^{-\nu} \Bigl( (-\nu) r^{-\nu-1} J_\nu(j_\nu r) 
+ r^{-\nu} j_\nu\, J_\nu'(j_\nu r) \Bigr) .
$$
Since $J_\nu (j_\nu) =0$ we obtain
\begin{equation} \label{e:phi'K}
\varphi_1'(1) \,=\, \kappa_n\, j_\nu^{-\nu+1}\, J_\nu'(j_\nu) .
\end{equation}
Furthermore,
$$
\varphi_1''(r) \,=\, \kappa_\nu \,j_\nu^{-\nu}
\Bigl( (-\nu) (-\nu-1) r^{-\nu-2} J_\nu(j_\nu r) +2(-\nu) r^{-\nu-1} j_\nu\, J_\nu'(j_\nu r) + r^{-\nu} j_\nu^2\, J_\nu'' (j_\nu r) \Bigr) 
$$
and hence
\begin{equation*} 
\varphi_1''(1) \,=\, \kappa_n \,j_\nu^{-\nu+1} \Bigl( -2\nu \, J_\nu'(j_\nu) + j_\nu \, J_\nu''(j_\nu) \Bigr).
\end{equation*}
To rewrite this further note that, by~\eqref{e:J2},
$$
2 J_\nu''(s) \,=\, J_{\nu-1}'(s) - J_{\nu+1}'(s) .
$$
Together with \eqref{e:J4} and \eqref{e:J3} we find
\begin{eqnarray*}
2s\,J_\nu''(s) &=&  s J_{\nu-1}'(s) - s J_{\nu+1}'(s) \\
&=&
\Bigl( (\nu-1) J_{\nu-1}(s) - s J_\nu (s) \Bigr) -
\Bigl( -(\nu+1) J_{\nu+1}(s) + s J_\nu (s) \Bigr) .
\end{eqnarray*}
At $s= j_\nu$ we obtain, together with \eqref{e:J1} and \eqref{e:J2},
$$
2 j_\nu\, J_\nu''(j_\nu) \,=\, J_{\nu+1}(j_\nu) - J_{\nu-1}(j_\nu) 
\,=\, -2 J_\nu'(j_\nu) .
$$
Altogether,
\begin{equation} \label{e:phi2fin}
\varphi_1''(1) \;=\; -\kappa_n\, j_\nu^{-\nu+1}\, (2\nu+1) J_\nu'(j_\nu) .
\end{equation}
In view of \eqref{e:sl}, \eqref{e:phi'K} and \eqref{e:phi2fin} 
the function $\gs_{\lef}(T)$ is equal to
\begin{equation} \label{sigmaleft}
\gs_{\lef} (T) \;=\; -\kappa_n\,j_\nu^{-\nu+1} J_\nu'(j_\nu) 
\Bigl( (2\nu+1) + \frac{ \xi I_{\nu+1}(\xi) }{ I_\nu(\xi) } \Bigr) .
\end{equation}
Using also \eqref{e:I3} and \eqref{e:I4} we can rewrite this as
\begin{equation} \label{sigmaleft2}
\gs_{\lef} (T) \;=\; -\kappa_n\,j_\nu^{-\nu+1} J_\nu'(j_\nu) 
\Bigl( 1 + \frac{ \xi I_{\nu-1}(\xi) }{ I_\nu(\xi) } \Bigr) .
\end{equation}
Since $\kappa_\nu$, $j_\nu$ are positive, $J_\nu'(j_\nu)$ is negative, 
and the functions $I_\nu$ are positive at all $\xi >0$,
formula~\eqref{sigmaleft2} implies

\begin{lemma} \label{le:sleft.g0}
In the interval of definition $(0,\mu)$ of the function $\gs_{\lef}$, we have
$$
\gs_{\lef} (T) >0 .
$$
\end{lemma}

Moreover, by~\eqref{e:power.I} we have
$$
\lim_{\xi \to 0} \frac{\xi I_{\nu+1}(\xi)}{I_\nu(\xi)} \,=\, 2 \, \frac{\Gamma (\nu+2)}{\Gamma (\nu+1)} .
$$
Since $\xi \to 0$ as $T \nearrow \mu$ by~\eqref{def:xi}, we find 
together with~\eqref{sigmaleft} that for all $\nu \ge 0$,
$$
\gs_1 (\mu) \,=\, \lim_{T \nearrow \mu} \gs_{\lef}(T) \,=\, 
-\kappa_n\,j_\nu^{-\nu+1} J_\nu'(j_\nu) 
\Bigl( 2 \nu + 1 + 2 \, \frac{\Gamma (\nu+2)}{\Gamma (\nu+1)} \Bigr) .
$$
In particular,
\begin{lemma} \label{le:satmu}
$\gs_1 (\mu) >0$.
\end{lemma}

\subsection{A formula for $\gs_{\righ}$.}
We follow the reasoning that we used to find a formula for 
the function $\gs_{\lef}(T)$. 
We skip the technical details. 
Assume that $T>\mu$. 
This allows us to define
\begin{equation} \label{def:rho}
\rho \,=\, \sqrt{\lambda_\nu - \left(\frac{2\,\pi}{T}\right)^2} \,.
\end{equation}
The function $\hat c (s) := s^\nu c \bigl( \frac s\rho \bigr)$ is 
the continuous solution on~$[0,\rho]$ of
$$
\left[\partial_s^2+ \frac{1}{s}\, \partial_s - \left(1 + \frac{\nu^2}{s^2}\right)\right]\,\hat c \,=\, 0
$$
with $\hat c(\rho) = -\rho^{\nu}\, \varphi_1'(1)$. 
The solution of this ordinary differential equation is given by 
$\beta\,J_{\nu}(s)$, where the constant $\beta$ (depending on $\nu$ and $T$)
is chosen such that 
$$
\beta\, J_\nu (\rho) \,=\, -\rho^{\nu}\, \varphi_1'(1) .
$$
Returning to the function~$c$, we get
$$
c(r) \,=\,  -\frac{\varphi_1'(1)}{J_\nu(\rho)}\, r^{-\nu}\, J_\nu(\rho\, r)
$$
and from \eqref{sigma:1:def} and \eqref{sigma:left:right}, using the identities~\eqref{e:J2}, 
\eqref{e:J3} and~\eqref{e:J4}, we obtain
\begin{eqnarray}
\gs_{\righ}(T) 
&=&
-\gf_1'(1) \frac{1}{J_\nu (\rho)} \frac 12 
\Bigl( \rho J_{\nu-1}(\rho) -2\nu J_\nu(\rho) - \rho J_{\nu+1}(\rho) \Bigr)  + \gf_1''(1) \notag \\
&=& 
\gf_1''(1) + \gf_1'(1) \frac{\rho J_{\nu+1}(\rho)}{J_\nu(\rho)} .  \label{e:sr}
\end{eqnarray}
In view of \eqref{e:phi'K} and \eqref{e:phi2fin} this becomes
\begin{eqnarray}  \label{e:gsRfin}
\gs_{\righ} (T) &=& -\kappa_n\,j_\nu^{-\nu+1} J_\nu'(j_\nu) 
\Bigl( (2\nu+1) - \frac{ \rho J_{\nu+1}(\rho) }{ J_\nu(\rho) } \Bigr) \\
&=& -\kappa_n\,j_\nu^{-\nu+1} J_\nu'(j_\nu) 
\Bigl( 1 + \frac{ \rho J_{\nu-1}(\rho) }{ J_\nu(\rho) } \Bigr) . \notag
\end{eqnarray}
where we used the identities \eqref{e:J3} and \eqref{e:J4} to get the second equality.


\section{Study of the derivative of $\gs_1(T)$}  \label{s:sigmaneg}

\ni
Throughout this section we assume again that $n \ge 2$. 
We start with

\begin{lemma} \label{lemma:asym}
The function $\gs_1 \colon (0,\infty) \to \RR$ has the asymptotics
$$
\lim_{T \to 0} \gs_1(T) = +\infty 
\qquad \text{ and } \qquad
\lim_{T \to \infty} \gs_1(T) = -\infty . 
$$
\end{lemma}

\proof
The first asymptotics is already proven in~\cite{Sicbaldi}.
We give an easier proof:
By~\eqref{def:xi} we have $\xi \to \infty$ as $T \to 0$.
Using~\eqref{sigmaleft2} and \eqref{asym:I} we therefore find 
$$
\lim_{T \to 0} \gs_1(T) \,=\, \lim_{\xi \to \infty} \frac{\xi \1 I_{\nu+1}(\xi)}{I_\nu(\xi)}
\,=\, \lim_{\xi \to \infty} \xi \,=\, \infty .
$$
To prove the second asymptotics, we read off from \eqref{def:rho} that
$\rho \nearrow \sqrt{\gl_\nu} = j_\nu$ as $T \to \infty$.
As is well-known, $j_\nu < j_{\nu+1}$ 
(see e.g.~\cite[\S 15$\cdot$22]{Watson}).
Therefore $J_{\nu+1}(j_\nu) >0$.
Together with~\eqref{e:gsRfin} we thus find
$$
\lim_{T \to \infty} \gs_1(T) \,=\, - \lim_{\rho \nearrow j_\nu} \frac{\rho J_{\nu+1}(\rho)}{J_\nu(\rho)}
\,=\, -\infty .
$$
as claimed.
\proofend

It is shown in~\cite[p.~336]{Sicbaldi} that the function $\gs_1$ is analytic and hence differentiable.
For our purposes, it would be enough to know that $\gs_1$ has exactly one zero $T_\nu$ and that
$\gs_1' (T_\nu) \neq 0$.
This follows from Lemma~\ref{le:sleft.g0}, Lemma~\ref{le:satmu}, Lemma~\ref{lemma:asym} and Lemma~\ref{l:right} below, 
that states that $\gs_1'(T) < 0$ for all $T \in (\mu,\infty)$.
We shall prove a somewhat stronger statement, namely that $\gs_1'(T) < 0$ for all $T \in (0,\infty)$.

\begin{proposition} \label{prop:sigmaneg}
Let $n \ge 2$. The function $\gs_1 \colon (0,\infty) \to \RR$ has negative derivative.
Moreover, $\gs_1$ has exactly one zero, say $T_\nu$.
\end{proposition}

\proof
We show that $\gs_{\lef}$ has negative derivative (Lemma~\ref{l:left}), 
that $\gs_{\righ}$ has negative derivative (Lemma~\ref{l:right}), 
and that $\sigma'_1(\mu) <0$ (Lemma~\ref{l:mu}). 
The fact that $\sigma_1$ has exactly one zero then follows together with Lemma~\ref{lemma:asym}.

\begin{lemma} \label{l:left}
$\gs'_{\lef}(T) <0$ for all $T \in (0,\mu)$.
\end{lemma}

\proof
Recall from \eqref{e:phi'K} that
$$
-\gf_1'(1) > 0 .
$$
Set $f(s) = \frac{sI_{\nu+1}(s)}{I_\nu(s)}$.
In view of~\eqref{e:sl} we need to show that $\frac{d}{dT} f \bigl( \xi(T)\bigr) <0$ for all $T \in (0,\mu)$.
Since $\frac{d}{dT} f (\xi(T)) = f'(\xi(T)) \,\xi'(T)$ and $\xi'(T)<0$ for all $T \in (0,\mu)$, 
this is equivalent to 
\begin{equation} \label{est:h}
f'(s) >0 \quad \mbox{ for all }\, s \in (0,\infty) .
\end{equation}
By \eqref{e:J3} we have 
$s I_{\nu+1}' = -(\nu+1) I_{\nu+1}+ sI_\nu$ and $sI_\nu' = -\nu I_\nu+sI_{\nu-1}$.
Therefore,
$$
f'(s) \,=\, \frac{s\bigl( I_\nu^2 - I_{\nu-1} I_{\nu+1} \bigr)}{I_\nu^2} .
$$
The lemma now follows from the following claim.

\begin{claim} \label{cl:In}
$I_\nu^2(s) > I_{\nu-1}(s)I_{\nu+1}(s)$ for all $\nu \in \RR$ and all $s > 0$.
\end{claim}

\proof
In view of \eqref{e:at0} we have $I_\nu^2(0) \ge I_{\nu-1}(0)I_{\nu+1}(0)$
for all $\nu \ge 0$. 
It therefore suffices to show that for all $s>0$,
$$
\frac{d}{ds} I_\nu^2 \,>\, \frac{d}{ds} \left( I_{\nu-1}I_{\nu+1}\right) .
$$
Multiplying by $s$, we see that this is true if and only if
\begin{equation} \label{e:2In}
2 I_\nu s I_\nu' \,>\, sI_{\nu-1}' I_{\nu+1} + I_{\nu-1}s I_{\nu+1}' .
\end{equation}
In view of \eqref{e:I2}, \eqref{e:I4}, \eqref{e:I3} we have
\begin{eqnarray*}
2sI_\nu' &=& sI_{\nu-1} + sI_{\nu+1} \\
sI_{\nu-1}' &=& \phantom{-}(\nu-1) I_{\nu-1} + sI_\nu \\
sI_{\nu+1}' &=& -(\nu+1) I_{\nu+1} + sI_\nu .
\end{eqnarray*} 
Therefore, \eqref{e:2In} holds if and only if 
$$
s I_{\nu-1}I_\nu + s I_\nu I_{\nu+1} \,>\, (\nu-1)I_{\nu-1}I_{\nu+1} + s I_\nu I_{\nu+1}
- (n+1) I_{\nu-1} I_{\nu+1} + s I_{\nu-1} I_\nu 
$$
i.e., 
$$
0 \,>\, -2 I_{\nu-1} I_{\nu+1}
$$
which is true because $I_\nu(s) >0$ for all $\nu \in \RR$ and $s>0$.
\proofend

\begin{lemma} \label{l:right}
$\gs'_{\righ}(T) <0$ for all $T \in (\mu,\infty)$.
\end{lemma}

\proof
Recall that $-\gf_1'(1) >0$.
Note that the function
$$
\rho \colon (\mu,\infty) \to (0,j_\nu), \quad
\rho(T) = \sqrt{\gl_\nu- \left( \frac{2\pi}{T} \right)^2}
$$
is strictly increasing.
Set $h(s) = \frac{sJ_{\nu+1}(s)}{J_\nu(s)}$.
In view of~\eqref{e:sr} we need to show that 
\begin{equation} \label{est:hJ}
h'(s) >0 \quad \mbox{ for all }\, s \in (0,j_\nu) .
\end{equation}
Since $j_\nu$ is the first positive zero of $J_\nu$, we see
as in the proof of Lemma~\ref{l:left} that~\eqref{est:hJ} is equivalent to
\begin{claim} \label{cl:Jn}
$J_\nu^2(s) > J_{\nu-1}(s)J_{\nu+1}(s)$ for all $s \in (0,j_\nu)$.
\end{claim}

\proof
Let again $j_{\nu-1}$, $j_{\nu}$, $j_{\nu+1}$ be the first positive zero of 
$J_{\nu-1}$, $J_{\nu}$, $J_{\nu+1}$, respectively.
Moreover, denote by $j_{\nu-1}^{(2)}$ the second positive zero of $J_{\nu-1}$.
Then
\begin{equation} \label{e:z}
j_{\nu-1} < j_\nu < j_{\nu+1}, \qquad j_\nu < j_{\nu-1}^{(2)} ,
\end{equation}
see e.g.~\cite[\S 15$\cdot$22]{Watson}.
It follows from the power series expansion~\eqref{exp:J} that
\begin{equation} \label{e:g0}
J_\nu(s) >0 \quad \mbox{ for }\, s \in (0,j_\nu) .
\end{equation}
Assume first that $s \in [j_{\nu-1},j_\nu)$.
Then~\eqref{e:z} and \eqref{e:g0} show that $J_\nu(s) >0$, $J_{\nu-1}(s) \le 0$, $J_{\nu+1}(s) >0$, whence the claim follows.
Assume now that $s \in (0,j_{\nu-1})$.
In view of~\eqref{eJ:at0} we have $J_\nu^2(0) \ge J_{\nu-1}(0) J_{\nu+1}(0)$
for all $\nu \ge 0$. 
It therefore suffices to show that 
\begin{equation}  \label{e:dds}
\frac{d}{ds} J_\nu^2 \,>\, \frac{d}{ds} \left( J_{\nu-1} J_{\nu+1}\right) 
\quad \mbox{ on }\, (0,j_{\nu-1}) .
\end{equation}
Using \eqref{e:J2}, \eqref{e:J4} and \eqref{e:J3}
we see as in the proof of Claim~\ref{cl:In} that~\eqref{e:dds} is equivalent to
$$
0 \,>\, -2 J_{\nu-1}(s)\, J_{\nu+1}(s) 
$$
which is true because $J_{\nu-1}$ and $J_{\nu+1}$ are positive on $(0,j_{\nu-1})$.
\proofend

To complete the proof of Proposition~\ref{prop:sigmaneg} we also show

\begin{lemma} \label{l:mu}
$\gs_1'(\mu) < 0$.
\end{lemma}

\proof
Since the function $\gs_1$ is smooth, 
$$
\gs_1'(\mu) \,=\, \lim_{T \searrow \mu} \gs'_{\righ} (T) .
$$
For $T>\mu$ we have $\gs_1'(T) = h'\bigl(\rho(T)\bigr) \, \rho'(T)$.
We compute
$$
\rho'(T) \,=\, \frac{2\pi}{\rho(T)\,T^3}
$$
and
\begin{equation} \label{e:h'}
h'(s) \,=\, \frac{s \left( J_\nu^2-J_{\nu-1}J_{\nu+1}\right)}{J_\nu^2} .
\end{equation}
Since $\lim_{T \searrow \mu}\rho(T) =0$ we obtain
$$
\gs_1'(\mu) \,=\, \lim_{T \searrow \mu} \gs'_{\righ}(T) \,=\, 
\frac{2\pi}{\mu^3} \,\gf_1'(1) \left( 1 - \lim_{s \to 0} \frac{J_{\nu-1} J_{\nu+1}}{J_\nu^2} \right) . 
$$
In view of the power series expansion~\eqref{exp:J},
$$
J_\nu(s) \,=\, \frac{(\frac 12 s)^\nu}{\Gamma(\nu+1)} + O (s^{s+\nu}) .
$$
Therefore, 
$$
\lim_{s \to 0} \frac{J_{\nu-1} J_{\nu+1}}{J_\nu^2} \,=\, \frac{\Gamma(\nu+1)}{\Gamma(\nu)\Gamma(\nu+2)} \,<\,1 \quad \mbox{ for all }\, \nu \ge 0
$$
and thus $\gs_1'(\mu) < 0$.
\proofend


\section{Extremal domains via the Crandall--Rabinowitz theorem} \label{s:extremalviaCR}

\ni
We are now in position to prove our main result when $n \geq 2$: 
The hypotheses of the Crandall--Rabinowitz bifurcation theorem are satisfied 
by the operator~$F$ defined in Section~\ref{s:proofI}. 
For $n \geq 2$, Theorem~\ref{t:main} follows at once from the following proposition and the Crandall--Rabinowitz theorem. 
As before, $\nu = \frac{n-2}{2}$. 

\begin{proposition} \label{cor:ker}
For $n \ge 2$, there exists a real number $T_*(n) = T_\nu$ 
such that the kernel of the linearized operator $D_v\,F(0,T_\nu)$ 
is 1-dimensional and is spanned by the function~$\cos t$,
$$
\Ker D_v\,F(0,T_\nu) \,=\, \RR \, \cos \1 t .
$$
The cokernel of $D_vF(0,T_\nu)$ is also 1-dimensional, and
$$
D_T\,D_v\,F(0,T_\nu) (\cos t) \,\notin\, \Im D_v\,F(0,T_\nu) .
$$
\end{proposition}

\proof
Let $v \in C^{2,\alpha}_{\textnormal{even},0}(\RR / 2 \pi \ZZ)$,
$$
v = \sum_{k \geq 1} a_k\, \cos(k\,t).
$$
We know that
\begin{equation}\label{aaa}
D_v\, F(0,T) \,=\, \sum_{k \geq 1} \sigma_k(T)\, a_k\, \cos(k\,t) .
\end{equation}
Let $V_k$ be the space spanned by the function $\cos(k\,t)$. 
By Proposition~\ref{prop:sigmaneg}, the function~$\sigma_1(T)$ has exactly one zero~$T_\nu$.  
By~\eqref{aaa}, the line~$V_1$ belongs to the kernel of $D_v\, F(0, T_\nu)$. 
Moreover, $V_1$ is the whole kernel, because for $k \ge 2$ we have
\[
\sigma_k(T_\nu) \,=\, \sigma_1\left(\frac{T_\nu}{k}\right) \neq 0
\]
(because $T_\nu$ is the only zero of $\sigma_1$). 
By~\eqref{aaa} 
and since $D_v\, F(0, T_\nu)$ is elliptic,
the image of $D_v\, F(0, T_\nu)$ is the closure of
\[
\bigoplus_{k\geq 2} V_k
\]
in $C^{1,\alpha}_{\textnormal{even},0}(\RR / 2 \pi \ZZ))$, and its codimension is equal to 1. 
More precisely,
\[
C^{1,\alpha}_{\textnormal{even},0}(\RR / 2 \pi \ZZ)) \,=\, \Im D_v\, F(0, T_\nu) \oplus V_1.
\]
Again by \eqref{aaa},
\[
D_T\, D_v\, F(0, T) (v) \,=\, \sum_{k \ge 1} \sigma_k'(T)\, a_k\, \cos(k\,t)
\]
and in particular
\[
D_T\, D_v\, F (0, T_\nu) (\cos t) \,=\, \sigma_1'(T_\nu)\, \cos t  \,\notin\,  \Im D_v\, F(0, T_\nu)
\]
because $\sigma_1'(T_\nu) < 0$ by Proposition~\ref{prop:sigmaneg}.
This completes the proof of the proposition. 
\proofend


\section{The problem in $\RR^2$}  \label{s:dim2}

\ni
Assume that $n=1$, i.e., the ambient space of the cylinder $C_1^T$ is $\RR^2$.
Recall from Section~\ref{s:proofI} that in this case,
\[
\gs_1(T) \,=\, c'(1) + \gf_1''(1)
\]
where $c$ is the solution of 
\begin{equation} \label{e:n=2}
\left( \pp_r^2 + \lambda_1 - \left(\frac{2\,\pi}{T}\right)^2\right)\, c \,=\, 0,
\end{equation}
with $c(1) = -\varphi_1'(1)$ and $c'(0)=0$,
where $\gf_1$ is the first eigenfunction of the Dirichlet problem on $[-1,1]$
normalized to have $L^2$-norm $\frac 1{2\pi}$.
(Here and in the sequel, $c$ denotes again the function $c_1$.)
For $\gf_1$ and $\gl_1$ we thus have 
\[
\gl_1 = \frac{\pi^2}{4} 
\qquad \textnormal{ and } \qquad 
\gf_1(r) = \frac{1}{\sqrt{2\pi}} \cos \left(\frac \pi 2 \2 r \right)
\] 
Hence
$$
-\gf_1'(1) = \sqrt{ \tfrac \pi 8} 
\quad \mbox{ and } \quad
\gf_1''(1) = 0 .
$$

\begin{lemma}
The only zero of the function $\gs_1 (T)$ is at $T = 4$. Moreover $\gs_1' (4) < 0$.
\end{lemma}

\proof
We abbreviate 
$\ga (T) := 
\gl_1-(\frac{2\pi}{T})^2 = \left( \frac \pi 2 \right)^2 - \left( \frac{2\pi}{T} \right)^2$.
The solution to~\eqref{e:n=2} is
$$
c(r) \,=\,
\left\{
\begin{array} {ll}
\sqrt{\frac \pi 8} \;  \frac{ \cosh \sqrt{-\ga (T)} \2 r }{ \cosh \sqrt{-\ga (T)} } 
  &\mbox{if }\; T \in (0,4) ,\\
\sqrt{\frac \pi 8}                     
  &\mbox{if }\; T = 4 ,\\
\sqrt{\frac \pi 8} \; \frac{ \cos \sqrt{\ga (T)} \2 r }{ \cos \sqrt{\ga (T)} }  
  &\mbox{if }\; T \in (4,\infty) .
\end{array} 
\right.
$$
Hence,
$$
\gs_1(T) = c'(1) \,=\,
\left\{
\begin{array} {ll}
-\sqrt{\frac \pi 8} \;  \sqrt{-\ga (T)} \tanh \sqrt{-\ga (T)}   
  &\mbox{if }\; T \in (0,4) ,\\
0                     
  &\mbox{if }\; T = 4 ,\\
-\sqrt{\frac \pi 8} \; \sqrt{\ga (T)} \tan \sqrt{\ga (T)} 
  &\mbox{if }\; T \in (4,\infty) .
\end{array} 
\right.
$$
In particular, $\gs_1 (T) >0$ on $(0,4)$ and $\gs_1(T) < 0$ on $(4,\infty)$.
It remains to show that $\gs_1' (4) < 0$.

For $T>4$ define $h(T) := \sqrt{\ga(T)}$. 
Then 
$$
\gs_1'(T) \,=\,
- \sqrt{\tfrac \pi 8}\, \tfrac{d}{dT} \bigl( h(T) \tan h(T) \bigr) \,=\,
- \sqrt{\tfrac \pi 8}\, h'(T) \bigl( \tan h(T) + h(T) ( 1+\tan^2 (h(T)) ) \bigr) .
$$
Since $\gs_1 (T)$ is smooth on $(0,\infty)$ and since $\displaystyle \lim_{T \to 4^+} h(T) =0$ and 
$h'(T) = \frac{\ga'(T)}{2 h(T)}$,
we find
\begin{eqnarray*}
\gs_1'(4) 
&=& 
- \sqrt{\tfrac \pi 8}\, \lim_{T \to 4^+} h'(T) \bigl( \tan h(T) + h(T) \bigr) \\
&=& 
- \sqrt{\tfrac \pi 8}\, \lim_{T \to 4^+} h'(T) \2 2 \2 h(T) \\
&=& 
- \sqrt{\tfrac \pi 8}\, \lim_{T \to 4^+} \ga'(T) 
\,=\,
- \sqrt{\tfrac \pi 8}\, \tfrac{\pi^2}8 
\,<\, 0 .
\end{eqnarray*}
\proofend

\begin{remark}
{\rm
A computation shows that $\gs_1' (T) < 0$ for all $T \in (0,\infty)$.
}
\end{remark}

Using the previous lemma, the proof of Proposition~\ref{cor:ker} applies also for $n=1$, 
and we obtain

\begin{proposition}\label{cor:ker:2}
Proposition~\ref{cor:ker} is true also for $n=1$ and $T_*(1) = 4$.
\end{proposition}

\begin{figure}[ht]
 \begin{center}
  \psfrag{t}{$t \in \RR$}
  \psfrag{x}{$x \in \RR$}
  \psfrag{-4}{$\approx -4$}
  \psfrag{4}{$\approx 4$}
 \leavevmode\epsfbox{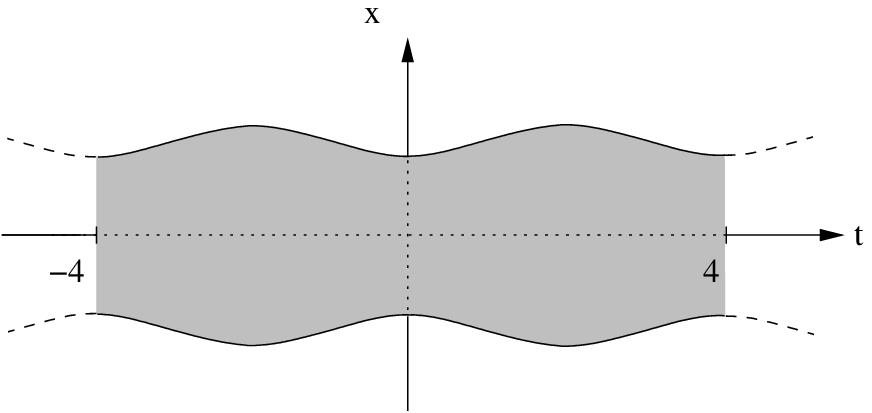}
 \end{center}
 \caption{A domain $\Omega_s \subset \RR^2$.}
 \label{figure.cyldim2}
\end{figure}

Together with the Crandall--Rabinowitz theorem we now obtain our main Theorem~\ref{t:main} 
also for $n=1$.
Figure~\ref{figure.cyldim2} shows the shape of the new extremal domains in~$\RR^2$.

\section{Estimates on the bifurcation period}  \label{s:pT}

\ni
Recall from Section~\ref{s:dim2} that $T_*(1) =4$.
In this section we study the bifurcation values $T_\nu = T_*(n)$ for $n \ge 2$,
and in particular prove Theorem~\ref{p:T}.

We recall that $J_\nu'(j_\nu) \neq 0$, and from~\eqref{est:hJ} 
that the function $h(s) = \frac{s J_{\nu+1}(s)}{J_\nu(s)}$ is strictly increasing on $(0,j_\nu)$ from~$0$ to~$\infty$.
By~\eqref{e:gsRfin} the unique zero $T_\nu$ of $\gs_{\righ}$ is therefore determined by
\begin{equation}  \label{e:rhonu}
\rho_\nu :=\ \rho (T_\nu) \,=\, \sqrt{\gl_\nu- \left(\frac{2\pi}{T_\nu}\right)^2}
\end{equation}
and
$$
\frac{\rho_\nu J_{\nu+1}(\rho_\nu)}{J_\nu(\rho_\nu)} \,=\, 2\nu+1 .
$$
In other words, the bifurcation value is 
\begin{equation} \label{e:bif}
T_\nu \,=\, \frac{2\pi}{\sqrt{\gl_\nu-\rho_\nu^2}}
\end{equation}
where $\rho_\nu$ is the unique zero on $(0,j_\nu)$ of 
$s J_{\nu+1} - (2\nu+1){J_\nu}$
or, by~\eqref{e:J1}, of $s J_{\nu-1} + J_\nu$.

For fixed $\nu$, the value $\rho_\nu$ and hence $T_\nu$ can be computed 
by the computer (using, for instance, Mathematica).
The first few and some larger values of~$T_\nu$ 
(rounded to five decimal places) are

\begin{equation} \label{t:10}
\begin{array}{|c||c|c|c|c|c|c|c|c|} 
\hline
2\nu         & 0             &1     & 2     &     3&     4&     5& 6    & 7    \\  \hline 
T_\nu &  3.06362 & 2.61931 & 2.34104 & 2.14351&1.99308 &1.87315 & 1.77429& 1.69088 \\
\hline 
2\nu         & 8             &9     & 10     &     11&     12&     13& 14    & 15    \\  \hline
T_\nu & 1.61924& 1.55650 & 1.50123 &1.45180 & 1.40735 &  1.36697& 1.33003&  1.2963\\
\hline
2\nu         & 16            & 17    & 18     & 19    & 20    & 40    & 200    & 2000    \\  \hline
T_\nu &  1.2650&   1.23616& 1.20927 & 1.18411&  1.16058 & 0.87348 &  0.4229& 0.13888  \\ 
\hline
\end{array}
\end{equation}

\b \ni
To study $T_\nu$ for $\nu \ge 10$ define 
$$
\rho_\nu^- = j_{\nu-1}+\frac{1}{j_{\nu-1}+2} , \qquad
\rho_\nu^+ = j_{\nu-1}+\frac{1}{j_{\nu-1}} .
$$

\begin{proposition} \label{p:T2}
The sequence $T_\nu$ is strictly decreasing to~$0$.
For $\nu \ge 10$ we have
\begin{equation} \label{e:est10}
\frac{2\pi}{\sqrt{\gl_\nu-(\rho_\nu^-)^2}} \,<\, T_\nu \,<\, \frac{2\pi}{\sqrt{\gl_\nu - (\rho_\nu^+)^2}} .
\end{equation}
\end{proposition}

\begin{remark} \label{rem:zeros}
{\rm
The zeros $j_\nu$ (and hence the eigenvalues $\gl_\nu = j_\nu^2$) 
are rather well-known, \cite{Lang-Wong}, namely
\begin{equation*} 
\nu - \frac{a_1}{\sqrt[3]{2}}\, \nu^{1/3} + \frac{3}{20} a_1^2\, \sqrt[3]{2}\, \nu^{-1/3} - 0.061\, \nu^{-1} \,<\, j_\nu \,<\, 
\nu - \frac{a_1}{\sqrt[3]{2}}\, \nu^{1/3} + \frac{3}{20} a_1^2\, \sqrt[3]{2}\, \nu^{-1/3} 
\end{equation*}
for all $\nu \in \frac 12 \NN$ with $\nu \ge 10$.
Here, $a_1 \approx -2.33811$ is the first negative zero of the Airy function $\Ai (x)$.
Therefore,
\begin{equation} \label{e:approx}
\nu + a\, \nu^{1/3} + b\, \nu^{-1/3} - c\, \nu^{-1} \,<\, j_\nu \,<\, 
\nu + a\, \nu^{1/3} + b\, \nu^{-1/3} 
\end{equation}
with positive constants $a \approx 1.8557$, $b \approx 1.0331$, $c < \frac 1{16}$. 
For $\gl_\nu$ we obtain the estimate
\begin{eqnarray} \label{e:approx.la}
\nu^2 + 2a\, \nu^{4/3} + \left( 2b+a^2 \right) \nu^{2/3} + 2ab + b^2\, \nu^{-2/3} - C(\nu) 
&<& \gl_\nu \;<\; \\
\nu^2 + 2a\, \nu^{4/3} + \left( 2b+a^2 \right) \nu^{2/3} + 2ab + b^2\, \nu^{-2/3} \phantom{-C(\nu) a} 
&& \notag
\end{eqnarray}  
where $C(\nu) = c \left( 2+2a\, \nu^{-2/3} +2b\, \nu^{-4/3} + c\, \nu^{-2} \right)$
is strictly decreasing, and $C(9) < 1/5$.
\diam
}
\end{remark}

\ni
We start with proving the estimate~\eqref{e:est10},
which by~\eqref{e:bif} is equivalent to
\begin{equation} \label{e:3rho}
\rho_\nu^- \,<\, \rho_\nu \,<\, \rho_\nu^+ .
\end{equation}
Recall that
\begin{equation}  \label{e:h2}
h(s) \,=\, \frac{sJ_{\nu+1}}{J_\nu} \,=\, 2\nu - \frac{sJ_{\nu-1}}{J_\nu} .
\end{equation}
Since $J_{\nu-1}(j_{\nu-1}) =0$
we have $h(j_{\nu-1}) = 2\nu$.
This and $h'(s) >0$ on $(0,j_\nu)$ show that
$$
j_{\nu-1} \,<\, \rho_\nu < j_\nu .
$$
In order to improve these bounds on $\rho_\nu$ we need
to better understand $h$ on the interval $I_\nu := [j_{\nu-1}, j_\nu)$.
The identities~\eqref{e:h'} and \eqref{e:h2} show that
\begin{equation}  \label{e:h'-2}
h'(s) \,=\, s+ \frac hs \left( h-2\nu \right) .
\end{equation}
In particular, 
\begin{equation}  \label{e:h'special}
h'(j_{\nu-1}) \,=\, j_{\nu-1} ,
\qquad
h'(\rho_\nu) \,=\, \rho_\nu + \frac{2\nu+1}{\rho_\nu} .
\end{equation}
Moreover, using~\eqref{e:h'-2}
\begin{eqnarray*}
h''(s) &=& 1+ \frac hs \,h' + \frac{h's-h}{s^2} (h-2\nu) \\
&=&
1+h+ (h-2\nu) \left( \frac{h^2-h}{s^2} +1 +\frac{h}{s^2} \,(h-2\nu) \right) .
\end{eqnarray*}
It follows that $h'' >0$ on $I_\nu$.
Therefore, the straight line of slope $h'(j_{\nu-1})$ passing through $(j_{\nu-1},2\nu)$
reaches the height $2\nu+1$ on the left of the graph of $h$,
while the straight line of slope $h'(\rho_\nu)$ passing through $(j_{\nu-1},2\nu)$
reaches the height $2\nu+1$ on the right of the graph of $h$, cf.~the figure below.

\begin{figure}[ht]
 \begin{center}
  \psfrag{j}{$j_{\nu-1}$}
  \psfrag{r-}{$>\rho_\nu^-$}
  \psfrag{r}{$\rho_\nu$}
  \psfrag{r+}{$\rho_\nu^+$}
  \psfrag{s}{$s$}
  \psfrag{2n}{$2\nu$}
  \psfrag{2n1}{$2\nu+1$}
  \psfrag{h}{$h(s)$}
 \leavevmode\epsfbox{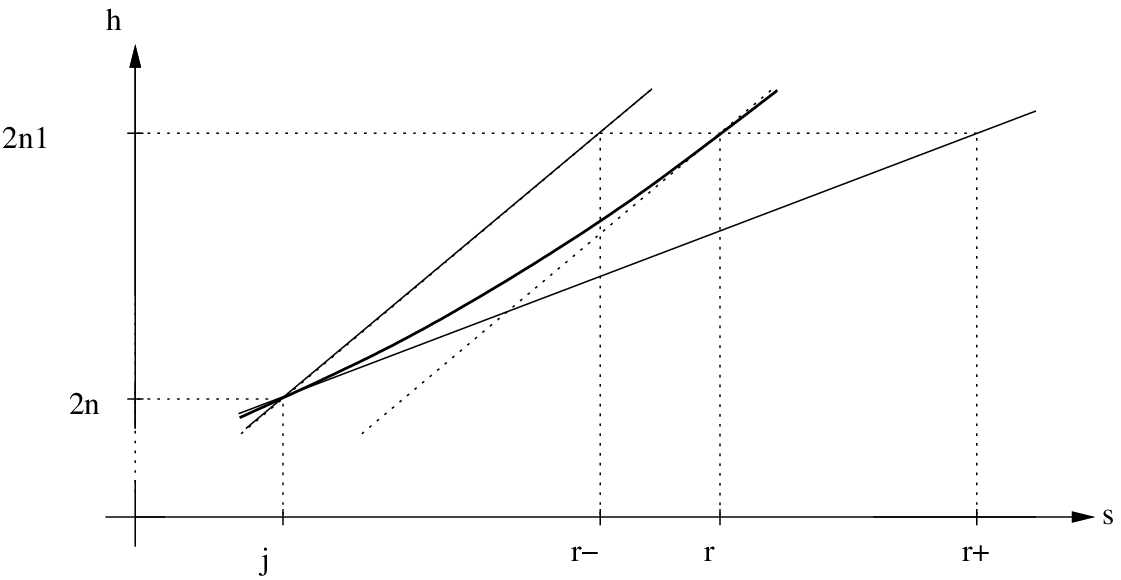}
 \end{center}
\end{figure}

\ni
Together with~\eqref{e:h'special} we conclude that 
$$
j_{\nu-1} + \frac 1{\rho_\nu+\frac{2\nu+1}{\rho_\nu}}  \,<\, \rho_\nu \,<\, j_{\nu-1} + \frac 1{j_{\nu-1}} =: \rho_\nu^+ .
$$
Using now that $j_{\nu-1} < \rho_\nu < j_{\nu-1} + \frac 1{j_{\nu-1}}$
and that $j_{\nu-1} > \nu+\frac 12$ for all $\nu \ge 10$ by~\eqref{e:approx},
we find
$$
\rho_\nu+\frac{2\nu+1}{\rho_\nu} \,<\, j_{\nu-1} + \frac{2\nu+2}{j_{\nu-1}} \,<\, j_{\nu-1}+2,
$$
and hence~\eqref{e:3rho} follows.

\m
It has been shown in~\cite{Sicbaldi} that 
$T_\nu < \frac{\sqrt 2 \pi}{\sqrt \nu}$,
whence $T_\nu$ converges to~$0$.
This also follows from~\eqref{e:est10} and
\begin{equation}  \label{e:O}
\gl_{\nu}-(\rho_\nu^+)^2 \,=\, j_{\nu}^2-j_{\nu-1}^2-2-\frac 1{j_{\nu-1}^2} \,=\,
2\nu + O \bigl( \nu^{1/3} \bigr)
\end{equation}
where for the last identity we used~\eqref{e:approx.la}.
Note that~\eqref{e:O} and $\gl_{\nu}-(\rho_\nu^-)^2 = 2\nu + O \bigl( \nu^{1/3} \bigr)$ imply that
$$
\frac{2\pi}{T_\nu} \,=\, \sqrt{2}\2 \nu^{1/2} + O(\nu^{-1/6}) 
\quad \mbox{ or } \quad
T_\nu \,=\, \sqrt{2}\2 \pi \2 \nu^{-1/2} +O(\nu^{-7/6}) .
$$
We finally show that the sequence~$T_\nu$ is strictly decreasing.
In view of the Table~\eqref{t:10} we can assume that $\nu \ge 10$.
By \eqref{e:est10} we need to show that for each such~$\nu$,
$$
\gl_\nu - \bigl( \rho_\nu^- \bigr)^2 \,<\, \gl_{\nu+\frac 12} - \bigl( \rho_{\nu+\frac 12}^+ \bigr)^2 ,
$$
i.e.,
\begin{equation} \label{e:convex}
\gl_{\nu +\frac 12} - \gl_{\nu} \,>\,  
\gl_{\nu -\frac 12} - \gl_{\nu-1} 
+ \left( 2- 2 \frac{j_{\nu-1}}{j_{\nu-1}+2} \right) + \left( \frac 1{j_{\nu-\frac 12}^2} - \frac 1{(j_{\nu-1}+2)^2} \right) .
\end{equation}
The first bracket on the RHS is equal to
$$
\frac 4{j_{\nu-1}+2} 
\,\stackrel{\tiny\eqref{e:approx}}{<}\, \frac 4{\nu+2} 
\,\le\, \frac 13,
$$
and the second bracket is less than $\frac 1{100}$.
It therefore suffices to show that
\begin{equation} \label{e:convex2}
\gl_{\nu +\frac 12} - \gl_{\nu} \,>\,  
\gl_{\nu -\frac 12} - \gl_{\nu-1} 
+ \frac 13 + \frac 1{100} .
\end{equation}
The function $\nu \mapsto \nu^\ga$ is convex for $\ga = \frac 43$ and $\ga = -\frac 23$,
but concave for $\ga = \frac 23$.
At $\nu =10$ we have
$$
(a^2+2b) \Bigl( (\nu+ \tfrac 12)^{2/3} - \nu^{2/3} \Bigr) \,>\,
(a^2+2b) \Bigl( (\nu- \tfrac 12)^{2/3} - (\nu-1)^{2/3} \Bigr) -\frac 1{30} .
$$
Furthermore,
$$
\left(\nu+ \tfrac 12 \right)^2 - \nu^2  \,=\,
\left(\nu- \tfrac 12 \right)^2 - (\nu-1)^2  +1 ,
$$
and $C(\nu-1) \le C(9) < \frac 15$ for $\nu \ge 10$.
Since $\frac 13 + \frac 1{100} + \frac 1{30} + 2 \frac 15 < 1$, 
the estimate~\eqref{e:approx.la} now implies that~\eqref{e:convex2} holds true. 
\proofend

\begin{remark}
{\rm
It is known that the function $\nu \mapsto \gl_\nu$ is strictly convex 
on $(0,\infty)$,
see~\cite{Elbert-Laforgia}.
In particular,
$$
\gl_{\nu+\frac 12} - \gl_\nu \,>\, 
\gl_{\nu-\frac 12} - \gl_{\nu-1} .
$$
This is not quite enough to prove inequality~\eqref{e:convex}.
}
\end{remark}

%
%
%


\enddocument